\documentstyle[12pt,leqno]{article}
\def\Bbb#1{{\hbox{\bf #1}}}
\def\R{\Bbb R}
\newcommand{\bbx}{\hfill\rule{2 mm}{2 mm}}
\newcommand{\inter}[1]{{\rm int}\,#1}
\newcommand{\spann}[1]{{\rm span}\,#1}

\begin{document}
\setcounter{page}{0}
\newcounter{aux}
\setcounter{aux}{6}
\title{The Volume of the Intersection of a Convex Body with its Translates}

\author{M.\ Meyer\\ Equipe d'Analyse\\Universit\'e Paris \Roman{aux} \\
F--75252--Paris, Cedex 05\\ France \\ \\
S.\ Reisner\\Department of Mathematics\\and\\School of Education of the
Kibbutz Movement--Oranim\\University of Haifa\\Haifa, 31905\\Israel \\ \\
M.\ Schmuckenschl\"ager\thanks{supported by Erwin Schr\"odinger Stipendium
J0630}\\
\begin{minipage}[c]{55 mm}
\begin{center}
Institut f\"ur Mathematik\\Johannes Kepler Universit\"at
\\Linz\\Austria\\ \mbox{}
\end{center}
\end{minipage}
\hfill and \hfill
\begin{minipage}[c]{65 mm}
\begin{center}
Department of Theoretical Mathematics\\
Weizmann Institute of Science\\Rehovot, 76100\\
Israel
\end{center}
\end{minipage}
\date{}}
\maketitle
\thispagestyle{empty}
\pagebreak

\begin{abstract}
It is proved that for a symmetric convex body $K$ in $\R^n$, if for some
$\tau
> 0$, $|K \cap (x+\tau K)|$ depends on $\left\|x\right\|_K$ only, then $K$
is an
ellipsoid. As a part of the proof, smoothness properties of convolution
bodies
are studied.\end{abstract}

\section{Introduction and statement of the main result }
Let $K$ be a symmetric convex body in $\R^n$ (i.e. $K$ is a compact, convex
set in $\R^n$ with non-empty interior, {\em symmetric} means in this paper:
centrally symmetric about the origin--0, unless another center of symmetry
is explicitely specified). We denote the norm (Minkowski-functional) which
is defined by $K$ by $\left\|\cdot\right\|_K$ :
$\|x\|_K=\inf \{\lambda>0\,;\,x\in \lambda K\}$.

We solve the following problem: What can be said about a symmetric convex
body $K$ which has the property that \mbox{$F(x)=|K \cap (x+K)|$} as a
function
of $x \in \R^n$, depends on $\left\|x\right\|_K$ alone? ( $|A|$ denotes the
$k$-dimensional Lebesgue measure of a set $A$, where $k$ is the dimension
of
the minimal flat containing $A$ ).

Actually, one could formulate a more general problem, using
$\left\|x\right\|_L$
where $L$ is another symmetric convex body. It is, however, easy to check,
which we do in the sequel, that the property formulated above, with
$\left\|x\right\|_L$ instead of $\left\|x\right\|_K$ implies that
\mbox{$\left\|x\right\|_L=\alpha\left\|x\right\|_K$}, i.e. $L$ and $K$ are
homothetic with center--0.

We prove in this paper that the only convex bodies having the above
property are ellipsoids.

The motivation for dealing with this kind of problems was to check to what
extent the following conjecture of Thue-Poulsen \cite{bib:TP}, Kneser
\cite{bib:K} and Hadwiger \cite{bib:H} can be generalized. The conjecture
(generalized to $\R^n$ ) is the following: Let $\{B_1$,$\ldots$,$B_k\}$,
$\{B'_1$,$\ldots$,$B'_k\}$ be two collections of Euclidean balls of the
same radius in $\R^n$. Let $x_i$, $x'_i$ be the centers of $B_i$, $B'_i$
respectively and assume that \mbox{$|x_i-x_j| \leq |x'_i-x'_j|$} for all
$i$ and $j$ ($|\cdot|$ denotes the Euclidean norm ).Then
\[\left|\bigcup_{i=1}^kB_i\right|\leq \left|\bigcup_{i=1}^kB'_i\right|\]
(cf. \cite{bib:B}, \cite{bib:G}, \cite{bib:C-P} and \cite{bib:Go}
 for solutions in special cases and related \linebreak results ).

>From our result it follows that the above conjecture can be formulated only
 about
Euclidean balls (or, which is the same, translates of an ellipsoid,
together with
the distance function defined by this ellipsoid ).

Actually the exact formulation of the result is stronger then the
one stated above
(cf. Theorem~\ref{th:1} ).

There are two main tools in the proof: one is a volumic formula for the
Gauss--Kronecker curvature of the boundary of a convex body, versions of which
were proved by Leichtweiss \cite{bib:L} and Schmuckenschl\"ager \cite{bib:S-m}
and which is extended here. From this formula and the volume-of-intersection
property assumed here, one gets a differential equation which can be transformed
into an equation of Monge--Amp\`er type and which was discussed by Petty
\cite{bib:P}. Petty's result is the second tool in the proof.

In order to obtain the result using the equation mentioned above, we had
to prove
sufficient smoothness of $K$. For this sake we prove some results, which are of
independent interest, concerning smoothness properties of {\em convolution
bodies.\/} The above name was given in \cite{bib:S-m} but such bodies were
defined and used before (cf. \cite{bib:M}, \cite{bib:B-L} ).

Let $K$ be a convex body in $\R^n$ with $0 \in\inter K$ and $\tau > 0$. For
$0 < \delta < \min (1,\tau^n)$ we set \[K(\delta,\tau) =\{x \in \R^n;|K
\cap (x+ \tau K)| \geq \delta \}\]

If $\tau =1$ we set $K(\delta)=K(\delta,1)$. We call $K(\delta,\tau)$--the
{\em convolution body\/} of $K$ with parameters $\delta$ and $\tau$ (this
is a deviation from the definition in \linebreak
\cite{bib:S-m} ). By Brunn--Minkowski
theorem, $K(\delta,\tau)$ is a convex body, moreover, it is symmetric if
$K$ is.

The structure of the paper is as follows: In the next parts of this
introduction
we give some notations and formulate the main theorem--Theorem~\ref{th:1}.
Then we bring a simple observation. Section~2 is devoted to proving
smoothness
properties of convolution bodies. In Section~3 we prove an extended volumic
formula for the Gauss--Kronecker curvature of the boundary of a convex
body.
In Section~4 we complete the proof of Theorem~\ref{th:1} and conclude with
some remarks and problems.
 \newline

The notations we use are standard.  The standard scalar
product in $\R^n$ is denoted by $\langle \cdot,\cdot \rangle$, the
Euclidean
unit sphere in $\R^n$ is $S^{n-1}$. If $K$ is a convex body in $\R^n$,
$\inter K$
is the interior of $K$ and $\partial K$ its boundary. If $0 \in \inter K$
we denote by $K^\ast$ the polar body of $K$ with respect to $0$ :
\[K^\ast = \{x \in \R^n\,; \forall y \in K, \langle x,y \rangle \leq 1\}
.\]
For a set $A \subset \R^n$, $\spann A$ denotes the linear span of $A$ and
\[A^\perp = \{x \in \R^n\,; \forall y \in A, \langle x,y \rangle =0\} .\]
If $0 \not=x \in \R^n$ we denote by $P_x$ the orthogonal projection onto
$\{x\}^\perp$.

\newtheorem{th:1}{Theorem}[section]
\begin{th:1}
Let $K$ be a symmetric convex body in $\R^n$. If for some $\tau >0$,
$|K \cap (x+ \tau K)|$ depends on $\|x\|_K$ only, Then (and only then) $K$
is an ellipsoid.
\label{th:1}
\end{th:1}
For $K$ and $L$ two symmetric convex bodies in $\R^n$ and $\tau >0$,
let $F(x)=|K \cap (x+ \tau K)|$.

\newtheorem{obs}[th:1]{Observation}
\begin{obs}
Let $K$ and $L$ be symmetric convex bodies in $\R^n$.
If $F(x)$ depends only on $\|x\|_L$ then $L$ is homothetic to $K$.
\label{obs}
\end{obs}

{\bf Proof.} Clearly $\{x;F(x)>0\}=\{x;\|x\|_K<1+ \tau \}$. Suppose now
that $K$ and $L$ are not homothetic. Then there exist $x,\,y$ such that
$\|x\|_L=\|y\|_L=1$ and $\|x\|_K=a>b=\|y\|_K$. Let $x'=((1+\tau)/a)x$,
 $y'=((1+\tau)/a)y$. Then $\|x'\|_K=1+\tau >\|y'\|_K$. It follows that
 $F(x')=0<F(y')$. But this contradicts the fact that
 $\|x'\|_L=\|y'\|_L=(1+\tau)/a$.\bbx

\section{Smoothness of Convolution--Bodies}

The results of this section are needed for the proof of Theorem \ref{th:1}.
But some of them are interesting independently of it.

Analogous results concerning strict convexity, smoothness and twice
differentiability
of the boundary of so-called {\em floating bodies\/}, were proved in
\cite{bib:M-R-2}.

\newtheorem{lem:2}{Lemma}[section]
\begin{lem:2}
Let $K_1$ and $K_2$ be two convex bodies in $\R^n$ and $u \in S^{n-1}$. For
$y \in P_uK_i$ $i=1,\,2$, denote: \pagebreak
\[\phi_i^+(y)=\max\{t;tu+y \in K_i\} \]
\[\phi_i^-(y)=\min\{t;tu+y \in K_i\} \]
\[f(r)=|K_1 \cap (ru+K_2)|\]
Then we have:
\[f'_+(0)=|{C_u}^+(1,2)|-|{C_u}^-(2,1)|\]
\[f'_-(0)=|{C_u}^-(1,2)|-|{C_u}^+(2,1)|\]
where
\[{C_u}^+(1,2)=P_u(K_1 \cap K_2) \cap \{\phi_1^+>\phi_2^+ \geq \phi_1^->
\phi_2^- \} \]
\[{C_u}^-(1,2)=P_u(K_1 \cap K_2) \cap \{\phi_1^+ \geq \phi_2^+>\phi_1^-
\geq \phi_2^-\} \]
and ${C_u}^\pm (2,1)$ are defined analogously.
\label{lem:2}
\end{lem:2}

{\bf Proof.} For $y \in P_u(K_1 \cap K_2)$ compute the limit at $y$, as
$r$ tends to $0$ from above, of
\begin{eqnarray*}
(1/r)\{[\min(\phi_1^+,\phi_2^++r)&-&\max(\phi_1^-,\phi_2^-+r)]_+\\
&-&[\min(\phi_1^+,\phi_2^+)-\max(\phi_1^-,\phi_2^-)]_+\}
\end{eqnarray*}
and integrate using Lebesgue theorem.\bbx
\newline

We bring now a proposition and a lemma. For the purpose of proving
Theorem~\ref{th:1}
Lemma~\ref{lem:4} would suffice. Proposition~\ref{prop:3} supplies,
however,
interesting information about convolution bodies in the case $\tau=1$,
which is not true in general in the case $\tau \not=1$.

\newtheorem{prop:3}[lem:2]{Proposition}
\begin{prop:3}
Let $K$ be a symmetric convex body in $\R^n$ and $0< \delta <|K|$. Then the
convolution body $K(\delta )$ is strictly convex.
\label{prop:3}
\end{prop:3}

{\bf Proof.} Assume $x_0 \not= x_1$, $\:x_0,\,x_1 \in K$, $x_t=x_0+t(x_1-
x_0)$ and $F(2x_t)=const.$ for $0 \leq t \leq 1$ ( $F(x)=|K \cap (x+K)|$ ).

Since $\{K_t=K \cap (2x_t+K)\}$ is a concave family, we conclude by the
equality
case in Brunn--Minkowski inequality and by symmetry with respect to $x_t$
of $K_t$, that \[K_t-x_t=K_0-x_0\]
Hence we get with $v=x_1-x_0$
\[K \supset K_0+[0,v]\;\;{\rm and}\;\;K \supset K_1+[0,v]\]
Therefore
\[[(K_0+[0,\;\;v]) \setminus K_0] \cap \inter{(2x_0+K)}=\emptyset\] and
\[[(K_1+[0,-v]) \setminus K_1] \cap \inter{(2x_1+K)}=\emptyset\]
It is easy to show that this implies
\[[(K_0+{\R_0}^+v) \setminus  K_0] \cap \inter{(2x_0+K)}=\emptyset\]
\[[(K_1+{\R_0}^+v) \setminus K_1] \cap \inter{(2x_1+K)}=\emptyset\]
Since $K_0$ is symmetric with respect to $x_0$ we have
\[2x_0+K \supset -K_0-[0,v]+2x_0=K_0-[0,v]\] and therefore
\[[(K_0-{\R_0}^+v) \setminus K_0] \cap \inter K=\emptyset\] and
\[[(K_1-{\R_0}^+v) \setminus K_1] \cap \inter K=\emptyset\]
These equations imply
\[(K_0+{\R}v) \cap K=K_0+[0,v]\;\;{\rm and}\;\;(K_0+{\R}v) \cap (2x_0+K)=
K_0-[0,v]\] Thus
\begin{eqnarray*}
(K_0+{\R}v) \cap (2x_t+K)&=&(K_0+{\R}v) \cap (K+2x_0+2tv)\\&=&K_0-[0,v]+2tv
\end{eqnarray*}
Hence \[K_{1/2}=K \cap (2x_{1/2}+K) \cap (K_0+{\R}v)=K_0+[0,v]\]
which is impossible.\bbx

\newtheorem{lem:4}[lem:2]{Lemma}
\begin{lem:4}
Let $K$ be a symmetric convex body such that for some $\tau >0$, $|K \cap
(x+\tau K)|$ depends on $\|x\|_K$ only. Then $K$ is strictly convex.
\label{lem:4}
\end{lem:4}

{\bf Proof.} We may assume $\tau <1$. Assume that $K$ is not strictly
convex, then there are points $x_0,\,x_1
\in \partial K$, $x_0$ an extreme point of $K$, such that for $0 \leq
\alpha
\leq 1$, $\|x_\alpha \|_K=1$ where $x_\alpha =(1-\alpha)x_0+\alpha x_1$.
By our
assumption $|K \cap (x_\alpha+\tau K)|$ is constant. Hence, by the equality
case in Brunn--Minkowski inequality, $K \cap (x_\alpha+\tau K)$ is a
translate of
$K \cap (x_0+\tau K)$. We show that \mbox{$K \cap (x_\alpha+\tau K)=\alpha
u+K \cap (x_0+\tau K)$} for $\alpha$ small enough, where \linebreak
$u=x_1-x_0$.

The point $(1-\tau)x_0$ is an interior point of $K$ and is an extreme point
 of $K \cap (x_0+\tau K)$. We can find an exposed point $y$ of $K \cap
 (x_0+\tau K)$, close enough to $(1-\tau)x_0$ so that $y \in \inter K$.

Let $v$ be an outer unit normal vector of $K \cap (x_0+\tau K)$ at $y$,
such
that $v$ is not an outer normal of $K \cap (x_0+\tau K)$ at any other
point.
Now, if $\alpha$ is small enough, $y+\alpha u \in \inter K$, therefore it
is an exposed point of $K \cap (x_\alpha+\tau K)$. $v$ is an outer normal
of
$K \cap (x_\alpha+\tau K)$ at $y+\alpha u$ and not at any other point. This
showes that it is impossible that \mbox{$K \cap (x_\alpha+\tau K)=w+K \cap
(x_0+\tau K)$} unless $w=\alpha u$ (since $v$ is an outer normal of $w+K
\cap (x_0+\tau K)$ at $y+w$).

Having shown this, it is sufficient to prove the lemma for the case $n=2$.
Suppose the boundary of $K$ is given locally by $\varphi :[-\varepsilon,
\varepsilon]\longrightarrow \R$ and
\[\varphi (t)=1\;\;{\rm for}\;\;t\geq 0\;,\;\;\;\;
\varphi (t)<1\;\;{\rm for}\;\;t<0\]
Obviously for $s>0$ the equality \mbox{$\varphi (t-s)=\varphi (t)$} is
impossible.\bbx

\newtheorem{th:5}[lem:2]{Theorem}
\begin{th:5}
Let $K$ be a convex body in $\R^n$ and $\tau >0$. If $K$ is strictly convex
then
for all $0< \delta < \min (1,\tau^n)|K|$, $\partial K(\delta,\tau)$ is of
class
${\cal C}^1$ (formula {\rm (\ref{eq:1})} in the following proof gives the
partial
derivative of $F(x)=|K \cap (x+\tau K)|$).
\label{th:5}
\end{th:5}

{\bf Proof.} By the implicit function theorem we have to show that
$F(x)$ has continuous partial derivatives on $\{x;F(x)>0\}$.

Let  $y \in \partial K$, $z \in \partial (x+\tau K)$, by $N(y)$ and $M(z)$
we
denote, respectively, the outer unit normal of $K$ at $y$ and the one of
$x+\tau K$
at $z$. By $\nu$ and $\mu$ we denote, respectively, the surface measures on
$\partial K$ and $\partial (x+\tau K)$. Note that $N$ and $M$ are well
defined a.e.\ $\nu$ (resp.\ $\mu$).
Since $K$ is strictly convex, $\partial K \cap \partial (x+\tau K)$ has
zero $(n-1)$-dimensional Hausdorff measure, hence \mbox{$|{C_u}^+(i,j)|=
|{C_u}^-(i,j)|$} in the notations of Lemma \ref{lem:2}. Therefore, for
\linebreak
$x=(x_1,\ldots,x_n)$ Lemma \ref{lem:2} can be written in the form
\begin{eqnarray}
\frac{\partial}{\partial x_1} F(x)&=&\int_{K \cap \partial (x+\tau K)}
\langle M(y),e_1
\rangle \,d\mu (y)\\ &=&-\int_{(x+\tau K) \cap \partial K}\langle N(y),e_1
\rangle
\,d\nu (y) \nonumber
\label{eq:1}
\end{eqnarray}
where $\{e_j\}$ are the coordinate vectors. The right hand side equality is
obtained by noticing that \[\frac{\partial}{\partial x_1}|K \cap (x+\tau K|
=-\frac{\partial}{\partial x_1}|\tau K \cap (-x+K)|\;.\]

Discontinuity may arise in (\ref{eq:1}) only in cases where for some $x$,
\linebreak
\mbox{$\partial
K \cap \partial (x+\tau K)$} has positive $(n-1)$-dimensional Hausdorff
measure,
by strict convexity, this does not happen, hence $\partial F/ \partial x_1$
is continuous on \linebreak
$\{x;F(x)>0\}$.

>From (\ref{eq:1}) it follows that at $x \in \partial K(\delta,\tau)$, i.e.\
such that $F(x)=\delta$, the outer normal of $\partial K(\delta,\tau)$ is
the well defined unit vector in the direction of
\begin{equation}
\int_{\partial K \cap (x+\tau K)} N(y)\,d\nu (y)
\label{eq:2}
\end{equation}
which is easily checked to be a non-zero vector for $\delta$ in the
specified range.\bbx

\newtheorem{th:6}[lem:2]{Theorem}
\begin{th:6}
Let $K$ be a convex body in $\R^n$ and $\tau >0$. If $K$ is strictly convex
and $\partial K$ is of class ${\cal C}^1$ then for $0<\delta <\min
(1,\tau^n) |K|$, $\partial K(\delta,\tau)$ is of class ${\cal C}^2$
 (formulae {\rm (\ref{eq:4})} and {\rm (\ref{eq:7})} in the following proof
 give the  second partial derivatives of $F(x)=|K \cap (x+\tau K)|$).
\label{th:6}
\end{th:6}

{\bf Proof.} Let us use the notations introduced in the proof of Theorem~
(\ref{th:5}). We compute first $\partial^2F/\partial x_2 \partial x_1$ for
$x \not= 0$ with $F(x)>0$. For $\varepsilon>0$ we have
\begin{equation}
\frac{1}{\varepsilon}\left[\frac{\partial}{\partial x_1}F(x+\varepsilon
e_2)-\frac{\partial}{\partial x_1}F(x)\right]\,=
\label{eq:3}
\end{equation}
\[-\frac{1}{\varepsilon}\left[\int_{\partial K \cap (x+\varepsilon e_2+
\tau K)} \langle N(y),e_1\rangle \,d\nu (y)\,-
\int_{\partial K \cap (x+\tau K)} \langle N(y),e_1 \rangle \,d\nu (y)
\right]\,=\]
\[=\frac{1}{\varepsilon}\int_{\partial K \cap [(x+\tau K)\bigtriangleup
(x+\varepsilon e_2+\tau K)]}\theta_\varepsilon (y) \langle N(y),e_1 \rangle
\,d\nu (y) \]
where \[\theta_\varepsilon(y)=-1\;\;{\rm if}\;\;y\in (x+\varepsilon e_2+
\tau K)\setminus (x+\tau K)\]
\[\theta_\varepsilon(y)=\;\;1\;\;{\rm if}\;\;y\in (x+\tau K)\setminus
(x+\varepsilon e_2+\tau K)\]

We claim that the limit as $\varepsilon$ tends to $0$ of the last integral
is
\begin{equation}
\frac{\partial^2}{\partial x_2\partial x_1}F(x)=-\int_{\partial K \cap
\partial (x+\tau K)}\frac{\langle M(y),e_2\rangle\langle N(y),e_1\rangle}
{[1\,-\,\langle M(y),N(y)\rangle^2]^{1/2}}\;d\sigma (y)
\label{eq:4}
\end{equation}
where $\sigma$ is the $(n-2)$-dimensional surface measure on \mbox
{$S=\partial K \cap \partial (x+\tau K)$} (which, as remarked in the proof
of Theorem \ref{th:5} has zero $(n-1)$-dimensional measure). We notice that
by the ${\cal C}^1$ assumption on $\partial K$, $N(y)$ and $M(y)$ are now
well defined and continuous on $S$.
By Stokes theorem or by a `reverse translation' argument like the one in
the proof of Theorem \ref{th:5}, we see that we can exchange the roles of
$M(y)$ and $N(y)$ in (\ref{eq:4}) and then get an expression which is
symmetric in $x_1$ and $x_2$ (see (\ref{eq:7})).

As $K$ is strictly convex, at each point $y$ of $S=\partial K \cap \partial
(x+\tau K)$, $N(y)$ and $M(y)$ are not parallel. It follows that \mbox
{$[1\,-\,\langle M(y),N(y)\rangle^2]^{1/2}$} is bounded away from $0$. At
every point $y\in S$ we estimate the length $L(y)$ of the $1$-dimensional
curve $\partial K \cap [y+E(y)]$, where \mbox{$E(y)=\spann{\{M(y),N(y)\}}
$}, in its intersection with the set $(x+\tau K) \bigtriangleup
(x+\varepsilon e_2+\tau K)$. We claim:
\begin{equation}
L(y)=\frac{\varepsilon |\langle M(y),e_2\rangle |}{[1\,-\,\langle M(y),N(y)
\rangle^2]^{1/2}}+\varepsilon \varphi_y(\varepsilon)
\label{eq:5}
\end{equation}
where $\varphi_y(\varepsilon)\longrightarrow 0$ as $\varepsilon
\longrightarrow 0$. This will be proved later.

For $\eta>0$ let \mbox{$S(\eta)=\{y\in S;|\langle M(y),e_2\rangle|\geq
\eta\}$} By compactness, for fixed $\eta$ and for small enough
$\varepsilon$ we have for all $y \in S(\eta)$
\begin{equation}
\theta_\varepsilon(z)=-sign~\langle M(y),e_2\rangle
\label{eq:6}
\end{equation}
for all $z\in\partial K \cap [(x+\tau K)\bigtriangleup(x+\varepsilon e_2+
\tau K)]\cap E(y)$.

The estimate (\ref{eq:5}) shows that the part of the last integral in
(\ref{eq:3}),
taken on the part of $\partial K \cap [(x+\tau K)\bigtriangleup(x+
\varepsilon
e_2+\tau K)]$ where $|\langle M(y),e_2\rangle|<\eta$ tends to $0$ as $\eta
\longrightarrow 0$. therefore, we conclude from (\ref{eq:6}), (\ref{eq:5})
and
(\ref{eq:3}), as $\varepsilon \longrightarrow 0$, the formula (\ref{eq:4}).

The same argument works for differentiation twice with respect to $x_1$, so
the formula for the second derivatives of $F$ may be written as
\pagebreak
\begin{equation}
\frac{\partial^2}{\partial x_i \partial x_j}F(x)=
\label{eq:7}
\end{equation}
\[\mbox{}-\frac{1}{2}\int_{\partial K \cap \partial(x+\tau K)}
\frac{\langle M(y),e_j\rangle\langle N(y),e_i\rangle+
\langle M(y),e_i\rangle\langle N(y),e_j\rangle}{[1\;-\;\langle M(y),N(y)
\rangle ^2]^{1/2}}\,d\sigma(y)\]

$i,\,j=1,\ldots,n$

As the integrand in (\ref{eq:7}) is continuous and $\partial K \cap
\partial
(x+\tau K)$ is a compact smooth $(n-2)$-dimensional manifold which varies
continuously, together with its derivatives on the surface $\partial K$,
all this by smoothness and strict convexity of $\partial K$, we get
continuity of the second derivatives.

Finally, to prove the estimate (\ref{eq:5}) we proceed as follows. Let
$\overline{n}$
be a unit vector in $\{N\}^\perp \cap \spann{\{M,N\}}$ (we drop the
argument $y$ here), and let $\overline{m}$ be a unit vector in $\{M\}^\perp
\cap \spann{\{\overline{n},e_2\}}$.

The approximation of first degree in $\varepsilon$ of $L(y)$ is the length
$l$
of the side parallel to $\overline{n}$ of the triangle which is formed by
lines parallel to $\overline{n},\,\overline{m}$ and $e_2$ and has an edge
of
length $\varepsilon$ in the direction of $e_2$. The estimate (\ref{eq:5})
follows now from an elementary (though somewhat lengthy) trigonometric
computation.
the degenerate case, where $\overline{m}$ and $\overline{n}$ are parallel
to $e_2$
is similarly shown to give $L(y)=o(\varepsilon)$ (this is of course the
case $\langle M,e_2\rangle=0$ ).\bbx
\newline

{\bf Remark concerning positive curvature.} In \cite{bib:M-R-2} it was
proved
that the floating body associated with a symmetric convex body $K$, where
$K$
is strictly convex and $\partial K$ is of class ${\cal C}^1$, is of class
${\cal C}^2$ and has positive Gauss--Kronecker  curvature at all points
of its boundary (this is not stated explicitely there, but the fact that
the
quadratic form $I(H,\cdot,\cdot)$ is positive-definite, proved in Corollary
2,
means this). It is probable that a similar fact is true for convolution
bodies.
Using Theorems \ref{th:5} and \ref{th:6} and in particular the formulae
(\ref{eq:2})
and (\ref{eq:7}) we can prove the following result, the details of whose
proof we omit.
\newtheorem{prop:7}[lem:2]{Proposition}
\begin{prop:7}
Let $K$ be a convex body in $\R^n$. If $\partial K$ is of class
${\cal C}^2$ and
has positive Gauss--Kronecker curvature everywhere, then the boundary of
$K(\delta)$ has
positive Gauss--Kronecker  curvature everywhere, for all $0<\delta<|K|$.
\label{prop:7}
\end{prop:7}

\section{A volumic formula for curvature}
\setcounter{equation}{0}

In this section we bring a formula for the generalized Gauss--Kronecker
curvature of the boundary of a convex body at a point where it exists. we
work
in a wider generality then what is needed for the proof of Theorem \ref
{th:1}
(for which one can assume that $\partial K$ is of class ${\cal C}^2$).

The terminology and notations of this section may be found e.g.\ in \cite
{bib:S-n} (pp.25--26) or in \cite{bib:L}. In particular, a {\em normal
point\/}
of the boundary of a convex body in $\R^n$, is a point where the Dupin
indicatrix
of $K$ exists. At a normal point, the Gauss--Kronecker curvature of $K$
exists and it can be
defined as the product of the principal curvatures (the reciprocal squares
of the semi-axes of the indicatrix) of $\partial K$ at the point $x$.

The set of points of $\partial K$ which are not normal is of zero $(n-1)$-
dimensional Hausdorff measure.

Let $K$ be a symmetric convex body in $\R^n$, $\tau>0$ and $x\in \partial
K$.
assume that the outer normal $N(x)$ of $K$ is well defined at $x$. For
$h>0$
small enough, denote by $x_h$ the point $\lambda x$ ($\lambda>0$) such that
$K \cap (x_h+\tau K)$ has width $h$ in the direction $N(x)$. The following
theorem is true also for non-symmetric $K$ if we replace $K \cap (x_h+\tau
K)$ by $K \cap (x_h-\tau K)$ and define $x_h$ accordingly.

\newtheorem{th:8}{Theorem}[section]
\begin{th:8}
Let $K$ be a symmetric convex body in $\R^n$ and $x$ a normal point of
$\partial
K$. For every $\tau>0$ the following formula gives the Gauss--\linebreak
Kronecker curvature $\kappa(x)=\kappa_K(x)$
\begin{equation}
\kappa(x)=\lim_{h \rightarrow 0} \frac{{c_n}^{n+1}h^{n+1}}{\left(\frac{1+
\tau}{\tau}\right)^{n-1}|K \cap (x_h+\tau K)|^2}
\label{eq:3.1}
\end{equation}
where $c_n=2\left(\frac{\omega_{n-1}}{n+1}\right)^\frac{2}{n+1}$ and
$\omega_k$ is the volume of the $k$-dimensional Euclidean ball.
\label{th:8}
\end{th:8}

{\bf Remarks.} {\em a)\/} It is proved in \cite{bib:L} (pp.450--454) that
at a normal point $x$
\begin{equation}
\kappa(x)=\lim_{h \rightarrow 0}\frac{{c_n}^{n+1}h^{n+1}}
{|K \cap{H_h}^+|^2}
\label{eq:3.2}
\end{equation}
where $H_h$ is a hyperplane orthogonal to $N(x)$, ${H_h}^+$-the half space
containing $x$ and bounded by $H_h$ and the width of $K \cap {H_h}^+$ in
direction $N(x)$ is $h$.

{\em b)\/} For $\tau=1$ (\ref{eq:3.1}) is a result of \cite{bib:S-m}.
\newline

{\bf Proof.} {\em case a):\/} $\kappa(x)>0$. The indicatrix of $K$ at $x$
is
an ellipsoid $E$. We choose an orthogonal coordinate system $x_1,\ldots,
x_n$
in $\R^n$ in which $\R^{n-1}=\{e_n\}^\perp$ is parallel to the tangent
hyperplane of $K$ at $x$. Using a volume-preserving linear transformation
which preserves $\R^{n-1}$, we may assume that $N(x)$ is parallel to $x$
and
that \[E=\left\{y\in \R^{n-1};\kappa(x)^\frac{1}{n-1}\sum_{j=1}^{n-1}
{y_j}^2 \leq 1\right\}\]
The indicatrix of $\tau K$ at $-\tau x$ is $\sqrt{\tau}E$.

Let $H_h$ and $J_h$ be hyperplanes parallel to $\R^{n-1}$ at distances $h$
and $\frac{\tau}{1+\tau}h$, respectively from $x$ (then $J_h$ is also of
distance $\frac{1}{1+\tau}h$ from $x_h-\tau x$ ).

By the definition of the indicatrix, the following limit exists in the
\linebreak Hausdorff-distance sense, as $h\longrightarrow 0$
\[\frac{1}{\sqrt{2h}} P_{e_n}(K \cap H_h)\longrightarrow E\]
>From which follows
\begin{equation}
\frac{1}{\sqrt{2h}} P_{e_n}(J_h \cap K) \longrightarrow \sqrt{\frac{\tau}
{1+\tau}} E
\label{eq:3.3}
\end{equation}
and by symmetry
\begin{equation}
\frac{1}{\sqrt{2h}} P_{e_n}(J_h \cap (x_h+\tau K)) \longrightarrow
\sqrt{\frac{\tau} {1+\tau}}E
\label{eq:3.4}
\end{equation}
Denote
\begin{eqnarray*}
V(h)&=&|K \cap (x_h+\tau K)| \\ V_1(h)&=&|K \cap {J_h}^+| \\
V_2(h)&=&|(x_h+\tau K) \cap {J_h}^-|
\end{eqnarray*}

>From (\ref{eq:3.3}) and (\ref{eq:3.4}) we conclude
\[|V(h)-[V_1(h)+V_2(h)]|\leq h^\frac{n+1}{2}f(h)\]
where $f(h) \longrightarrow 0$ as $h \longrightarrow 0$. We get
\begin{equation}
\lim_{h \rightarrow 0} \frac{V(h)}{{c_n}^\frac{n+1}{2}h^\frac{n+1}{2}}=
\lim_{h \rightarrow 0} \frac{V_1(h)+V_2(h)}{{c_n}^\frac{n+1}{2}h^\frac{n+1}
{2}}
\label{eq:3.5}
\end{equation}

Using (\ref{eq:3.5}), the distance from $J_h$ to $x$ and to $x_h-\tau x$,
(\ref{eq:3.2}),
and the identity $\kappa_{\tau K}(-\tau x)=\frac{1}{\tau}\kappa(x)$ we get
\[\lim_{h \rightarrow 0} \frac{V(h)}{{c_n}^\frac{n+1}{2}h^\frac{n+1}{2}}=
\frac{\left(\frac{\tau}{1+\tau}\right)^\frac{n+1}{2}}{\sqrt{\kappa(x)}}+
\frac{\left(\frac{1}{1+\tau}\right)^\frac{n+1}{2}\tau^\frac{n-1}{2}}
{\sqrt{\kappa(x)}}=
\left(\frac{\tau}{1+\tau}\right)^\frac{n-1}{2}\kappa(x)^{-1/2}
\] which proves (\ref{eq:3.1}) in {\em case a)\/}.

{\em case b):\/} $\kappa(x)=0$. After an appropriate transformation, we may
assume
that the indicatrix $E$ of $K$ at $x$ is of the form \mbox{$E=E_{n-1-j}
\times \R^j$} where
\[E_{n-1-j}=\left\{y \in \R^{n-1-j};
k^\frac{1}{n-1-j}\sum_{i=1}^{n-1-j}{y_i}^2 \leq 1\right\}\]

For $\beta>0$ let
\[S(h)=\{y \in \R^{n-1};-\beta\sqrt{2h}\leq y_i \leq \beta\sqrt{2h}\;\;{\rm
for}\;\;i=n-j,\ldots,n-1\}\]

Let $J_h$ be as in {\em case a)\/} and
\begin{eqnarray*}V(h)=|K \cap (x_h+\tau K)|\;&\hspace{5 mm}&\;W(h)=|K \cap
(x_h+\tau K)\cap S(h)| \\
W_1(h)=|K \cap S(h) \cap {J_h}^+|&\hspace{5 mm}&W_2(h)=|(x_h+\tau K)\cap
S(h) \cap {J_h}^-|
\end{eqnarray*}

By the definition of the indicatrix, $(1/\sqrt{2h})P_{e_n}(K \cap J_h)$ and
\linebreak
$(1/\sqrt{2h})P_{e_n}((x_h+\tau K) \cap J_h)$ tend to infinity in the
directions
$(y_{n-j},\ldots,y_{n-1})$, as $h \longrightarrow 0$ and this limit is
uniform
in the subspace of $\R^{n-1}$ orthogonal to these directions. Also, $S(h)$
cuts
the above projection sets, in these directions at $\pm \beta$ (when $h$ is
small
enough). Therefore we get, like in {\em case a)\/}:
\[|W(h)-[W_1(h)+W_2(h)]| \leq h^\frac{n+1}{2}\cdot o(1)\;\;{\rm as}\;\;h
\longrightarrow 0.\] hence
\begin{eqnarray}
\liminf_{h \rightarrow 0}\frac{V(h)}{{c_n}^\frac{n+1}{2}h^\frac{n+1}{2}}
&\geq&
\liminf_{h \rightarrow 0}\frac{W(h)}{{c_n}^\frac{n+1}{2}h^\frac{n+1}{2}}\\
&=&\liminf_{h \rightarrow 0}\frac{W_1(h)+W_2(h)}{{c_n}^\frac{n+1}{2}h^
\frac{n+1}{2}} \nonumber
\label{eq:3.6}
\end{eqnarray}

As $h \longrightarrow 0\;$ $W_1(h)$ is, up to terms of the form $o(h^
\frac{n+1}{2})$, not less then
\[\beta^j|E_{n-1-j}|\int_0^\frac{\tau}{1+\tau}(2t)^\frac{n-1}{2}\,dt=
\gamma h^\frac{n+1}{2}\beta^j\] where $\gamma$ does not depend on $h$ or
$\beta\;\;$
(cf. the computation in \cite{bib:L} p.452). With similar estimate for
$W_2(h)$ we get
\[\liminf_{h \rightarrow 0}\frac{V(h)}{{c_n}^\frac{n+1}{2}h^\frac{n+1}{2}}
\geq
\gamma'\beta^j\] Since this is true for all $\beta>0$ {\em case b)\/} is
proved.\bbx

\section{Proof of Theorem \protect\ref{th:1}}
\setcounter{equation}{0}
Let $K$ and $\tau$ satisfy the assumptions of Thoerem \ref{th:1}. Then for
$0<\delta< \min(1,\tau^n)|K|$, $K$ is homothetic to $K(\delta,\tau)$. By
Lemma \ref{lem:4} and Theorems \ref{th:5} and \ref{th:6}, it follows that
$K$ is strictly convex with boundary of class ${\cal C}^2$.

The Gauss map $u \longleftrightarrow x(u)$ between $u \in S^{n-1}$ and the
point
$x(u) \in \partial K$ where the outer normal of $K$ is $u$, is well defined
and
invertible. Using this map we can write the Gauss--Kronecker curvature of
$K$ as
a function $\kappa(u)$ defined on $S^{n-1}$. In this notation, Theorem
\ref{th:8} is written:
\begin{equation}
\kappa(u)=\lim_{\lambda \nearrow 1}\frac{{c_n}^{n+1}(1-\lambda)^{n+1}
|\langle x(u),u \rangle|^{n+1}}{\left(\frac{1+\tau}{\tau}\right)^{n-1}
|K \cap [(1+\lambda)x(u)+\tau K]|^2}
\label{eq:4.1}
\end{equation}

Clearly $\langle x(u),u \rangle=\|u\|_{K^\ast}$. By the assumption, and
since
$\|x(u)\|_K=1$, $|K \cap [(1+\lambda)x(u)+\tau K]|$ depends on $\lambda$
but not on $u$.

Let $f(u)=(1/ \kappa(u))=R_1\cdot \ldots \cdot R_{n-1}$ be the product of
the principal radii of curvature of $K$ at $x(u)$. We get
\begin{equation}
f(u)=\frac{k}{{\|u\|_{K^\ast}}^{n+1}}
\label{eq:4.2}
\end{equation}
for some constant $k$. By the Gauss--Bonnet formula, the constant $k$
cannot be
infinite hence $f(u)$ is finite everywhere and is a curvature function in
the
sense of being a derivative of the curvature measure of $K\;$
(cf. \cite{bib:S-n}).
\pagebreak

The following theorem was proved by Petty \cite{bib:P} (Lemma 8.4), using
deep results of Pogorelov, Cheng and Yau and others

\newtheorem{th:9}{Theorem}[section]
\begin{th:9}[Petty]
Let $K$ be a convex body in $\R^n$ with $\partial K$ of class ${\cal C}^2$
and having a curvature function $f(u)$. If $f(u)$ satisfies $(2)$
then $K$ is an ellipsoid.
\label{th:9}
\end{th:9}

Using Theorem \ref{th:9}, the proof of Theorem \ref{th:1} is concluded.\bbx
\newline

{\bf Concluding remarks.} {\em a)\/} We discussed here in the context of
Theorem \ref{th:1}, only centrally symmetric convex bodies. The proprty of
convex bodies which is assumed in Theorem \ref{th:1} can be formulated for
a-priory not necessarily symmetric bodies. One should specify then a point
$0 \in \inter K$ to be the origin. For the case $\tau=1$ it turns out that
convex
bodies which satisfy the assumption of Theorem \ref{th:1} must be
symmetric,
hence, by our result, ellipsoids. Actually, it follows from Lemma
(\ref{lem:2}) that such a convex body $K$ must be
homothetic to the polar of the projection body of $K\;$ (cf.\ also
\cite{bib:S-m} Theorem 1).

{\em b)\/} Remark {\em a)\/} above brings one back to the well known
problem
of Busemann: If for $n \geq 3$, a convex body $K$ in $\R^n$ is homothtic to
the
polar of its projection body, must $K$ be an ellipsoid? A positive answer
to this problem would give another solution to
 the problem which is solved in this paper (for
$ \tau=1$ and $n \geq 3$) by using the assumption $|K \cap (x+K)|=
f(\|x\|_K)$
only for small translations $x$. Our proof uses this assumption only for
`large'
translations, i.\ e.\ vectors $x$ such that $\|x\|_K$ is close to 2.

In fact, one can ask whether the following property characterizes
ellipsoids:
{\em For one fixed $0< \alpha <1+\tau$, $|K \cap (x+\tau K)|$ is constant
for $\|x\|_K=\alpha$\/}.

{\em c)\/} The case $\tau=\infty$, i.\ e.\ assumption on the volumes of
sections
of $K$ by hyperplanes in all directions can be proved along the same lines
of
the proof of Theorem \ref{th:1}, using floating bodies instead of
convolution
bodies. A very simple proof of this case is given however in
\cite{bib:M-R-1}
(Lemma 3) (though the proof here has the slight advantage of using the
assumption only `near the boundary' of $K$).

{\em d)\/} A weaker result concerning the problem treated in this paper had
been proved before by Y.~Gordon and M.~Meyer (unpublished). They proved
that if for {\em every} $\tau >0$, $|K \cap (x+\tau K|$ depends only on
$\|x\|_K$, then $K$ is an ellipsoid.
\pagebreak

\end{document}